\documentclass[a4paper,leqno]{article}

\usepackage{mathrsfs}
\usepackage[ansinew]{inputenc}
\usepackage{float}
\usepackage{afterpage}
\usepackage{slashbox}
\usepackage{amsmath,amssymb,amsthm,amsfonts}

\usepackage{amsmath,amssymb,amsfonts}
\usepackage{graphicx}

\def\RR{{{\mathbb R}}}
\def\ZZ{{{\mathbb Z}}}
\def\NN{{{\mathbb N}}}
\newtheorem{theorem}{Theorem}[section]

\title{A fractional spline collocation-Galerkin method for the time-fractional diffusion equation}
\author{Laura Pezza\thanks{ {\it Dept. SBAI, University of Roma ''La Sapienza''}, 
		Via A. Scarpa 16, 00161 Roma, Italy.
		e-mail: {\tt laura.pezza@sbai.uniroma1.it}},
	Francesca Pitolli\thanks{{\it Dept. SBAI, University of Roma ''La Sapienza''}, 
		Via A. Scarpa 16, 00161 Roma, Italy. 
		e-mail: \tt{francesca.pitolli@sbai.uniroma1.it}} 
}

\date{}

\begin{document}
	
\maketitle

\begin{abstract}
The aim of this paper is to numerically solve a diffusion differential problem having time derivative of fractional order.
To this end we propose a collocation-Galerkin method that uses the fractional splines as approximating functions. The main advantage is in that the derivatives of integer and fractional order  of the fractional splines can be expressed in a closed form that involves just the generalized finite difference operator. 
This allows us to construct an accurate and efficient numerical method. 
Several numerical tests showing the effectiveness of the proposed method are presented.
\end{abstract}

\section{Introduction.}
\label{sec:intro}

The use of fractional calculus to describe real-world phenomena is becoming increasingly 
widespread. Integro-differential equations of {\em fractional}, {\em i.e.} positive real, order
are used, for instance, to model wave propagation in porous materials, 
diffusive phenomena in biological tissue, viscoelastic properties of continuous media \cite{Hi00,Ma10,KST06,Ta10}.
Among the various fields in which fractional models are successfully used, viscoelasticity is one of the more
interesting since the memory effect introduced by the time-fractional  derivative allows to model
anomalous diffusion phenomena in materials that have mechanical properties in between pure elasticity and pure viscosity \cite{Ma10}.
Even if these models are empirical, nevertheless they are shown to be consistent with experimental data. 
\\
The increased interest in fractional models has led to the development of several numerical methods to solve
fractional integro-differential equations. Many of the proposed methods 
generalize to the fractional case numerical methods commonly used for the classical integer case
(see, for instance, \cite{Ba12,PD14,ZK14} and references therein). But the nonlocality
of the fractional derivative raises the challenge of obtaining numerical solution with high
accuracy at a low computational cost. In \cite{PP16} we proposed a collocation method
especially designed for solving differential equations of fractional order in time. 
The key ingredient of the method is the use of the fractional splines introduced in \cite{UB00} as approximating functions.
Thus, the method takes advantage of the explicit differentiation rule for fractional B-splines that allows us to evaluate accurately the derivatives 
of both integer and fractional order. 
\\
In the present paper we used the method to solve a diffusion problem having time derivative
of fractional order and show that the method is efficient and accurate.
More precisely, the {\em fractional spline collocation-Galerkin method} here proposed combines  the fractional spline collocation method introduced in \cite{PP16} for the time discretization and a classical spline Galerkin method 
in space.  
\\
The paper is organized as follows. 
In Section~\ref{sec:diffeq}, a time-fractional diffusion problem is presented and the definition of fractional derivative is given. 
Section~\ref{sec:fractBspline} is devoted to the fractional B-splines and the explicit expression of their fractional derivative is given.
The fractional spline approximating space is described in Section~\ref{sec:app_spaces},
while the fractional spline collocation-Galerkin method is introduced in Section~\ref{sec:Galerkin}. 
Finally, in Section~\ref{sec:numtest} some numerical tests showing the performance of the method are displayed.
Some conclusions are drawn in Section~\ref{sec:concl}.

\section{A time-fractional diffusion problem.}
\label{sec:diffeq}

We consider the {\em time-fractional differential diffusion problem} \cite{Ma10} 
\begin{equation} \label{eq:fracdiffeq}
\left \{ \begin{array}{lcc} 
\displaystyle D_t^\gamma \, u(t, x) - \frac{\partial^2}{\partial x^2} \, u(t, x) = f(t, x)\,, & \quad t \in [0, T]\,, &  \quad x \in [0,1] \,,\\ \\
u(0, x) = 0\,, &  & \quad x \in [0,1]\,, \\ \\
u(t, 0) = u(t, 1) = 0\,, & \quad t \in [0, T]\,,    
\end{array} \right.
\end{equation}
where $ D_t^\gamma u$, $0 < \gamma < 1$, denotes the {\em partial fractional derivative} with respect to the time $t$.
Usually, in viscoelasticity the fractional derivative is to be understood in the Caputo sense, {\em i.e.}
\begin{equation} \label{eq:Capfrac}
D_t^\gamma \, u(t, x) = \frac1{\Gamma(1-\gamma)} \,  \int_0^t \, \frac{u_t(\tau,x)}{(t - \tau)^\gamma} \, d\tau\,, \qquad t\ge 0\,,
\end{equation} 
where $\Gamma$ is the Euler's gamma function
\begin{equation}
\Gamma(\gamma+1)= \int_0^\infty \, s^\gamma \, {\rm e}^{-s} \, ds\,.
\end{equation}
We notice that due to the homogeneous initial condition for the function $u(t,x)$, solution of the differential problem (\ref{eq:fracdiffeq}), 
the Caputo definition (\ref{eq:Capfrac}) coincides with the Riemann-Liouville definition  (see \cite{Po99} for details). 
One of the advantage of the Riemann-Liouville definition is in that the usual differentiation operator in  the Fourier domain can be easily extended 
to the fractional case, {\em i.e.}
\begin{equation}
{\cal F} \bigl(D_t^\gamma \, f(t) \bigr) = (i\omega)^\gamma {\cal F} (f(t))\,,
\end{equation}
where ${\cal F}(f)$ denotes the Fourier transform of the function $f(t)$. Thus, analytical Fourier methods usually used in the classical integer case
can be extended to the fractional case \cite{Ma10}.

\section{The fractional B-splines and their fractional derivatives.}
\label{sec:fractBspline}

The {\em fractional B-splines}, {\em i.e.} the B-splines of fractional degree, were introduced in \cite{UB00} generalizing to the fractional power the classical definition
for the polynomial B-splines of integer degree.
Thus, the fractional B-spline $B_{\alpha}$ of degree $\alpha$   is defined as
\begin{equation} \label{eq:Balpha}
B_{\alpha}(t) := \frac{{ \Delta}^{\alpha+1} \, t_+^\alpha} {\Gamma(\alpha+1)}\,, \qquad \alpha > -\frac 12\,,
\end{equation}
where
\begin{equation} \label{eq:fracttruncpow}
t_+^\alpha: = \left \{ \begin{array}{ll}
t^\alpha\,, & \qquad t \ge 0\,, \\ \\
0\,, & \qquad \hbox{otherwise}\,,
\end{array} \right. \qquad \alpha > -1/2\,,
\end{equation}
is the {\em fractional truncated power function}. $\Delta^{\alpha}$ is the {\em generalized finite difference operator}
\begin{equation} \label{eq:fracfinitediff}
\Delta^{\alpha} \, f(t) :=  \sum_{k\in \NN} \, (-1)^k \, {\alpha \choose k} \, f(t-\,k)\,, 
\qquad \alpha \in \RR^+\,,
\end{equation}
where
\begin{equation}  \label{eq:binomfrac}
{\alpha \choose k} := \frac{\Gamma(\alpha+1)}{k!\, \Gamma(\alpha-k+1)}\,, \qquad k\in \NN\,, \quad \alpha \in \RR^+\,,
\end{equation}
are the {\em generalized binomial coefficients}. We notice that 'fractional' actually means 'noninteger', {\em i.e.} $\alpha$ can assume any real value 
greater than $-1/2$. For real values of $\alpha$, $B_\alpha$ does not have compact support even if it belongs to $L_2(\RR)$. 
When $\alpha=n$ is a nonnegative integer, Equations~(\ref{eq:Balpha})-(\ref{eq:binomfrac}) are still valid; $\Delta^{n}$ is 
the usual finite difference operator so that $B_n$ is the classical polynomial B-spline of degree $n$ and compact support $[0,n+1]$
(for details on polynomial B-splines see, for instance, the monograph \cite{Sc07}). 

\begin{figure}[t]
\begin{center}
\begin{tabular}{cc}
\includegraphics[width=6cm]{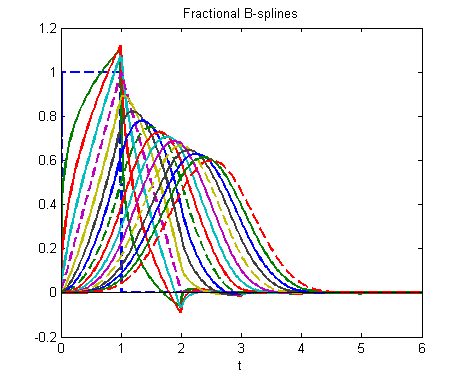}   & \includegraphics[width=6cm]{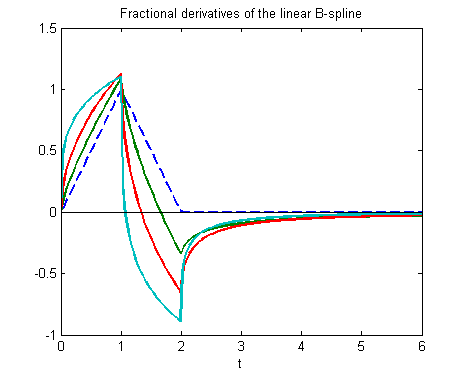} \\
\includegraphics[width=6cm]{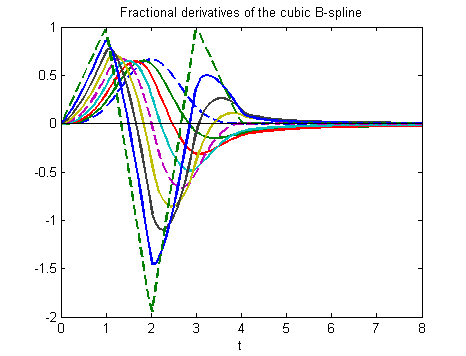} & \includegraphics[width=6cm]{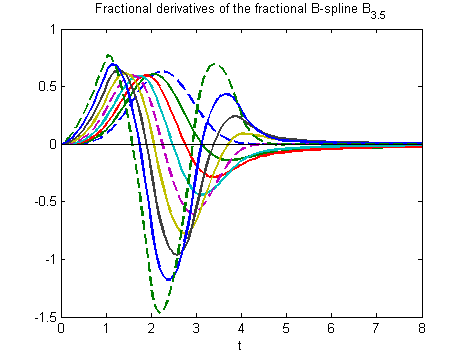}
\end{tabular}
\end{center}
\caption{Top left panel: The fractional B-splines (solid lines)  and the polynomial B-splines (dashed lines) for $\alpha$ ranging from 0 to 4. 
Top right panel: The fractional derivatives of the linear B-spline $B_1$ for $\gamma = 0.25, 0.5, 0.75$. 
Bottom left panel: The fractional derivatives of the cubic B-spline $B_3$ for $\gamma$ ranging from 0.25 to 2. Bottom right panel: The fractional derivatives of the fractional B-spline $B_{3.5}$ for the $\gamma$ ranging from 0.25 to 2. Ordinary derivatives are displayed as dashed lines. 
}
\label{fig:fractBsplines}
\end{figure}

The fractional B-splines for different values of the parameter $\alpha$ are displayed in Figure~\ref{fig:fractBsplines} (top left panel). 
The classical polynomial B-splines are also displayed (dashed lines). The picture shows that
the fractional B-splines decay very fast toward infinity so that they can be assumed compactly supported for computational 
purposes. Moreover, in contrast to the polynomial B-splines, fractional splines are not always positive even if the nonnegative part becomes more and more smaller as $\alpha$ increases. 

The fractional derivatives of the fractional B-splines can be evaluated explicitly by differentiating (\ref{eq:Balpha}) and (\ref{eq:fracttruncpow}) in the Caputo sense.
This gives the following differentiation rule 
\begin{equation}\label{eq:diffrule_tronc}
D^{\gamma}_t \, B_{\alpha} (t)=  \frac{\Delta^{\alpha+1} \, t_+^{\alpha-\gamma}} {\Gamma(\alpha-\gamma+1)}\,, \qquad 0 < \gamma < \alpha + \frac12\,,
\end{equation}
which holds both for fractional and integer order $\gamma$. In particular, when $\gamma, \alpha$ are nonnegative integers, (\ref{eq:diffrule_tronc})
is the usual differentiation rule for the classical polynomial B-splines \cite{Sc07}. 
We observe that since $B_\alpha$ is a causal function with $B_\alpha^{(n)}(0)=0$ for $n\in \NN\backslash\{0\}$, 
the Caputo fractional derivative coincides with the Riemann-Liouville fractional derivative.
\\
From (\ref{eq:diffrule_tronc}) and the composition property $\Delta^{\alpha_1} \, \Delta^{\alpha_2} = \Delta^{\alpha_1+\alpha_2}$
 it follows \cite{UB00}
\begin{equation} \label{eq:diffrule_2}
D^{\gamma}_t \, B_{\alpha} = \Delta ^{\gamma} \, B_{\alpha-\gamma}\,,
\end{equation}
{\em i.e.} the fractional derivative of a fractional B-spline of degree $\alpha$ is a fractional spline of degree $\alpha-\gamma$.
The fractional derivatives of the classical polynomial B-splines $B_n$ are fractional splines, too. 
This means that $D^{\gamma}_t \, B_{n}$ is not compactly supported when $\gamma$ is noninteger reflecting
the nonlocal behavior of the derivative operator of fractional order. 
\\
In Figure~\ref{fig:fractBsplines} the fractional derivatives of $B_1$ (top right panel), $B_3$ (bottom left panel) and $B_{3.5}$ (bottom right panel)
are displayed for different values of $\gamma$. 

\section{The fractional spline approximating spaces.}
\label{sec:app_spaces}

A property of the fractional B-splines that is useful for the construction of numerical methods for the solution of differential problems is the {\em refinability}. 
In fact, the fractional B-splines are {\em refinable functions}, {\em i.e.} they satisfy the {\em refinement equation}
\begin{equation}
B_\alpha(t) = \sum_{k\in \NN} \, a^{(\alpha)}_{k} \, B_\alpha(2\,t-k)\,, \qquad t \ge 0\,,
\end{equation}
where the coefficients
\begin{equation}
a^{(\alpha)}_{k} := \frac{1}{2^{\alpha}} {\alpha+1 \choose k}\,,\qquad k\in \NN\,,
\end{equation}
are the {\em mask coefficients}. This means that the sequence of nested approximating spaces
\begin{equation}
V^{(\alpha)}_j(\RR) = {\rm span} \,\{B_\alpha(2^j\, t -k), k \in \ZZ\}\,, \qquad j \in \ZZ\,,
\end{equation}
forms a {\em multiresolution analysis} of $L_2(\RR)$. 

As a consequence, any function $f_j(t)$ belonging to $V^{(\alpha)}_j(\RR)$ can be expressed as 
\begin{equation}
f_j(t) = \sum_{k\in \ZZ}\, \lambda_{jk} \, B_\alpha(2^j\, t -k)\,,
\end{equation}
where the coefficient sequence $\{\lambda_{j,k}\} \in \ell_2(\ZZ)$.  
Moreover, any space $V^{(\alpha)}_j(\RR)$ reproduces polynomials up to degree $\lceil \alpha\rceil$, 
{\em i.e.} $x^d \in V^{(\alpha)}_j(\RR)$, $ 0 \le d \le \lceil \alpha\rceil$,
while its approximation order is $\alpha +1$. We recall that the polynomial B-spline $B_n$ reproduces polynomial up to degree $n$ whit approximation order
$n+1$ \cite{UB00}.

To solve boundary differential problems we need to construct a multiresolution analysis on a finite interval. For the sake of simplicity in the following 
we will consider the interval $I=[0,1]$.  
A simple approach is to restrict the basis $\{B_\alpha(2^j\, t -k)\}$ to the interval $I$, {\em i.e.}
\begin{equation} \label{eq:Vj_int}
V^{(\alpha)}_j(I) = {\rm span} \,\{B_\alpha(2^j\, t -k), t\in I, -N \le k \le 2^j-1\}\,, \qquad j_0 \le j\,,
\end{equation}
where $N$ is a suitable index, chosen in order the significant part of $B_\alpha$ is contained in $[0,N+1]$, and $j_0$ is the starting refinement level.
The drawback of this approach is its numerical instability and the difficulty in fulfilling the boundary conditions since there are $2N$
boundary functions, {\em i.e.} the translates of $B_\alpha$ having indexes $ -N\le k \le -1$ and $2^j-N\le k \le 2^j-1$, that are non zero at the boundaries.
More suitable refinable bases can be obtained by the procedure given in \cite{GPP04,GP04}. In particular, for the polynomial B-spline $B_n$ a B-basis 
$\{\phi_{\alpha,j,k}(t)\}$ with optimal approximation properties can be constructed. 
The internal functions $\phi_{\alpha,j,k}(t)=B_\alpha(2^j\, t -k)$, $0 \le k \le 2^j-1-n$, remain unchanged
while the $2n$ boundary functions fulfill the boundary conditions
\begin{equation}
\begin{array}{llcc}
\phi_{\alpha,j,-1}^{(\nu)}(0) = 1\,, & \phi_{\alpha,j,k}^{(\nu)}(0) = 0\,, &\hbox{for} & 0\le \nu \le -k-2\,, -n \le k \le -2\,,\\ \\
\phi_{\alpha,j,2^j-1}^{(\nu)}(1) = 1\,, & \phi_{\alpha,j,2^j+k}^{(\nu)}(1) = 0\,, & \hbox{for} & 0\le \nu \le -k-2\,, -n \le k \le -2\,,\\ \\
\end{array}
\end{equation}
Thus, the B-basis naturally fulfills Dirichlet boundary conditions.
\\
As we will show in the next section, the refinability of the fractional spline bases plays a crucial role in the construction of the collocation-Galerkin method.

\section{The fractional spline collocation-Galerkin method.}
\label{sec:Galerkin}

In the collocation-Galerkin method here proposed, we look for an approximating function
$u_{s,j}(t,x) \in V^{(\beta)}_s([0,T]) \otimes  V^{(\alpha)}_j([0,1])$.
Since just the ordinary first spatial derivative of $u_{s,j}$ is involved in the Galerkin method, we can assume $\alpha$ integer 
and use as basis function for the space $V^{(\alpha)}_j([0,1])$ the refinable B-basis $\{\phi_{\alpha,j,k}\}$, {\em i.e.}
\begin{equation} \label{uj}
u_{s,j}(t,x) = \sum_{k \in {\cal Z}_j} \, c_{s,j,k}(t) \, \phi_{\alpha,j,k}(x)\,,
\end{equation}
where the unknown coefficients $c_{s,j,k}(t)$ belong to $V^{(\beta)}_s([0,T])$. Here, ${\cal Z}_j$ denotes the set of indexes $-n\le k \le 2^j-1$.
\\
The approximating function $u_{s,j}(t,x)$ solves the variational problem
\begin{equation} \label{varform}
\left \{ 
\begin{array}{ll} 
\displaystyle \left ( D_t^\gamma u_{s,j},\phi_{\alpha,j,k} \right ) -\left ( \frac {\partial^2} {\partial x^2}\,u_{s,j},\phi_{\alpha,j,k} \right ) = \left ( f,\phi_{\alpha,j,k} \right )\,, & \quad k \in {\cal Z}_j\,, \\ \\
u_{s,j}(0, x) = 0\,, & x \in [0,1]\,, \\ \\
u_{s,j}(t, 0) = 0\,, \quad u_{s,j}(t,1) = 0\,, & t \in [0,T]\,,
\end{array} \right.
\end{equation}
where $(f,g)= \int_0^1 \, f\,g$. 
\\
Now, writing (\ref{varform}) in a weak form and using (\ref{uj}) we get the system of fractional ordinary differential equations
\begin{equation} \label{fracODE}
\left \{ \begin{array}{ll}
M_j \, D_t^\gamma\,C_{s,j}(t) + L_j\, C_{s,j}(t) = F_j(t)\,, & \qquad t \in [0,T]\,, \\ \\
C_{s,j}(0) = 0\,,
\end{array} \right.
\end{equation}
where  $C_{s,j}(t)=(c_{s,j,k}(t))_{k\in {\cal Z}_j}$ is the unknown vector. 
The connecting coefficients, i.e. the entries of the mass matrix $M_j = (m_{j,k,i})_{k,i\in{\cal Z}_j}$, of the stiffness matrix $L_j = (\ell_{j,k,i})_{k,i\in{\cal Z}_j}$,
and of the load vector $F_j(t)=(f_{j,k}(t))_{k\in {\cal Z}_j}$,   are given by
$$
m_{j,k,i} = \int_0^1\, \phi_{\alpha,j,k}\, \phi_{\alpha,j,i}\,, \qquad \ell_{j,k,i} = \int_0^1 \, \phi'_{\alpha,j,k} \, \phi'_{\alpha,j,i}\,,
$$
$$ 
f_{j,k}(t) = \int_0^1\, f(t,\cdot)\, \phi_{\alpha,j,k}\,.
$$
The entries of $M_j$ and $L_j$ can be evaluated explicitly using (\ref{eq:Balpha}) and (\ref{eq:diffrule_tronc}), respectively, while the entries of $F_j(t)$ 
can be evaluated by quadrature formulas especially designed for wavelet methods \cite{CMP15,GGP00}.

\medskip
To solve the fractional differential system (\ref{fracODE}) we use the collocation method introduced in \cite{PP16}. 
For an integer value of $T$, let $t_p =  p/2^q$, $0\le p \le 2^q\,T$, 
where $q$ is a given nonnegative integer,
be a set of dyadic nodes in the interval $[0,T]$. Now, assuming 
\begin{equation} \label{ck}
c_{s,j,k}(t) = \sum_{r\in {\cal R}_s} \, \lambda_{k,r}\,\chi_{\beta,s,r}(t) \,, \qquad k \in {\cal Z}_j\,,
\end{equation} 
where $\chi_{\beta,s,r}(t)=B_\beta(2^s\,t-r)$ with $B_\beta$ a fractional B-spline of fractional degree $\beta$, 
and collocating (\ref{fracODE}) on the nodes $t_p$, we get the linear system
\begin{equation} \label{colllinearsys}
(M_j\otimes A_s + L_j\otimes G_s) \,\Lambda_{s,j}  =F_j\,, 
\end{equation}
where $\Lambda_{s,j}=(\lambda_{k,r})_{r\in {\cal R}_s,k\in {\cal Z}_j}$ is the unknown vector,
$$ 
\begin{array}{ll}
A_s= \bigl( a_{p,r} \bigr)_{p\in {\cal P}_q,r\in {\cal R}_s}\,, & \qquad a_{p,r} = D_t^\gamma \, \chi_{\beta,s,r}(t_p)\,, \\ \\
G_s=\bigl(g_{p,r}\bigr)_{p \in {\cal P}_q,r\in {\cal R}_s}\,, & \qquad g_{p,r} = \chi_{\beta,s,r}(t_p)\,, 
\end{array}
$$
are the collocation matrices and
$$
F_j=(f_{j,k}(t_p))_{k\in{\cal Z}_j,p \in {\cal P}_q}\,, 
$$
is the constant term. Here, ${\cal R}_s$ denotes the set of indexes $-\infty < r \le 2^s-1$ and ${\cal P}_q$ denotes the set of indexes $0<p\le 2^qT$.
Since the fractional B-splines have fast decay, 
the series (\ref{ck}) is well approximated by only few terms and 
the linear system (\ref{colllinearsys}) has in practice  finite dimension so that the unknown vector $\Lambda_{s,j}$ can be recovered by solving 
(\ref{colllinearsys}) in the least squares sense.  
\\
We notice that the entries of $G_s$, which involve just the values of $\chi_{\beta,s,r}$ on the dyadic nodes $t_p$, can be evaluated explicitly by (\ref{eq:Balpha}).
On the other hand, we must pay a special attention to the evaluation of the entries of $A_s$ since they involve the values of the fractional derivative 
$D_t^\gamma\chi_{\beta,s,r}(t_p)$. As shown in Section~\ref{sec:fractBspline}, they can be evaluated efficiently by the differentiation rule (\ref{eq:diffrule_2}).

In the following theorem we prove that the fractional spline collocation-Galerkin method is convergent.
First of all, let us introduce the Sobolev space on bounded interval
$$
H^\mu(I):= \{v \in L^2(I): \exists \, \tilde v \in H^\mu (\RR) \ \hbox{\rm such that} \ \tilde v|_I=v\}, \quad \mu\geq 0\,,
$$
equipped with the norm
$$
\|v\|_{\mu,I}  = \inf_{\tilde v \in H^\mu(\RR), \tilde v|_I=v} \|\tilde v\|_{\mu,\RR}\,,
$$
where
$$
H^\mu(\RR):= \{v: v\in L^2(\RR) \mbox{ and } (1+|\omega|^2)^{\mu/2}  {\cal F}(v)(\omega) \in L^2(\RR)\}, \quad \mu\geq 0\,,
$$
is the usual Sobolev space with the norm
$$
\| v \| _{\mu,\RR} =\bigl \| (1+|\omega|^2)^{\mu/2}  {\cal F}(v)(\omega) \bigr \| _{0,\RR}\,.
$$

\begin{theorem} \label{Convergence}
Let 
$$
H^\mu(I;H^{\tilde \mu}(\Omega)):= \{v(t,x): \| v(t,\cdot)\|_{H^{\tilde \mu}(\Omega)} \in H^\mu(I)\}, \quad \mu, {\tilde \mu} \geq 0\,,
$$
equipped with the norm
$$
\|v\|_{H^\mu(I;H^{\tilde \mu}(\Omega))} := \bigl \| \|v(t,\cdot)\|_{H^{\tilde \mu}(\Omega)} \bigr\|_{\mu,I}\,.
$$
Assume $u$ and $f$ in (\ref{eq:fracdiffeq}) belong to $H^{\mu}([0,T];H^{\tilde \mu}([0,1]))$, $0\le \mu$, $0\le \tilde \mu$, and $H^{\mu-\gamma}([0,T];$ $H^{{\tilde \mu}-2}([0,1]))$, $0\le \mu-\gamma$, $0\le \tilde \mu-2$, respectively. 
Then, the fractional spline collocation-Galerkin method is convergent, {\em i.e.},
\begin{equation}
\|u-u_{s,j}\|_{H^0([0,T];H^0([0,1]))} \, \to 0  \quad \hbox{as} \quad s,j \to \infty\,. 
\end{equation}
Moreover, for $\gamma \le \mu \le \beta+1$ and $1 \le \tilde \mu \le \alpha +1$ the following error estimate holds:
\begin{equation}
\begin{array}{lcl}
\| u-u_{s,j}\|_{H^0([0,T];H^0([0,1]))} &\leq & \left (\eta_1 \, 2^{-j\tilde \mu}  + \eta_2  \, 2^{-s\mu} \right ) \| u\|_{H^\mu([0,T];H^{\tilde \mu}([0,1]))}\,, 
\end{array}
\end{equation}
where $\eta_1$ and $\eta_2$ are two constants independent of $s$ and $j$.
\end{theorem}

\begin{proof}	
Let $u_j$ be the exact solution of the variational problem (\ref{varform}).  Following a classical line of reasoning (cf. \cite{Th06,FXY11,DPS94}) we get
$$
\begin{array}{l}
\| u-u_{j,s}\|_{H^0([0,T];H^0([0,1]))} \leq \\ \\ 
\rule{2cm}{0cm} \leq \|u-u_{j}\|_{H^0([0,T];H^0([0,1]))} + \| u_j-u_{j,s}\|_{H^0([0,T];H^0([0,1]))}\, \leq 
\\ \\ \rule{2cm}{0cm} \leq \eta_1  \, 2^{-j\tilde \mu}\, \| u\|_{H^0([0,T];H^{\tilde \mu}([0,1]))} + \eta_2 \, 2^{-s\mu} \, \| u\|_{H^\mu([0,T];H^0([0,1]))} \leq 
\\ \\ \rule{2cm}{0cm} \leq  \left ( \eta_1  \, 2^{-j\tilde \mu} + \eta_2 \, 2^{-s\mu} \right ) \, \|u\|_{H^\mu([0,T];H^{\tilde \mu}([0,1]))}\,.
\end{array}
$$
\end{proof}

\section{Numerical tests.}
\label{sec:numtest}

\begin{figure}[htp]
	\centering
	\begin{tabular}{cc}
		\includegraphics[width=0.45\textwidth]{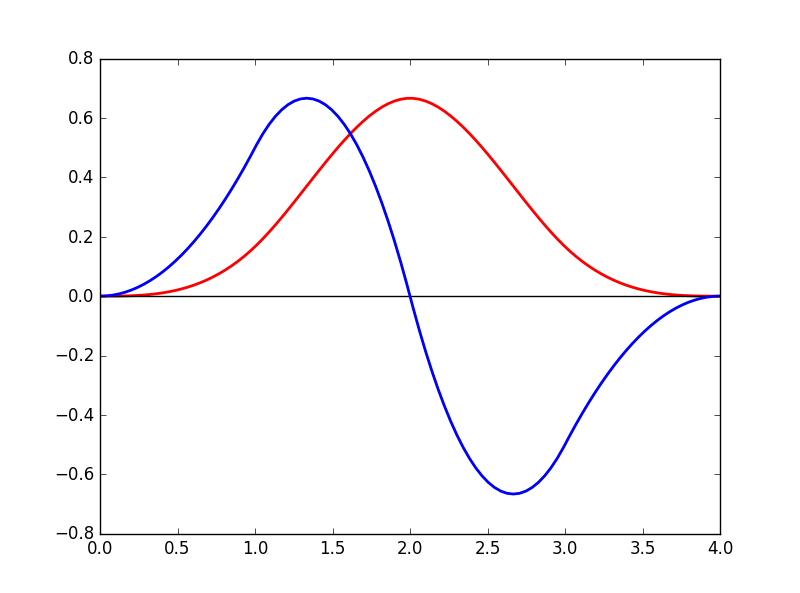} &
		\includegraphics[width=0.45\textwidth]{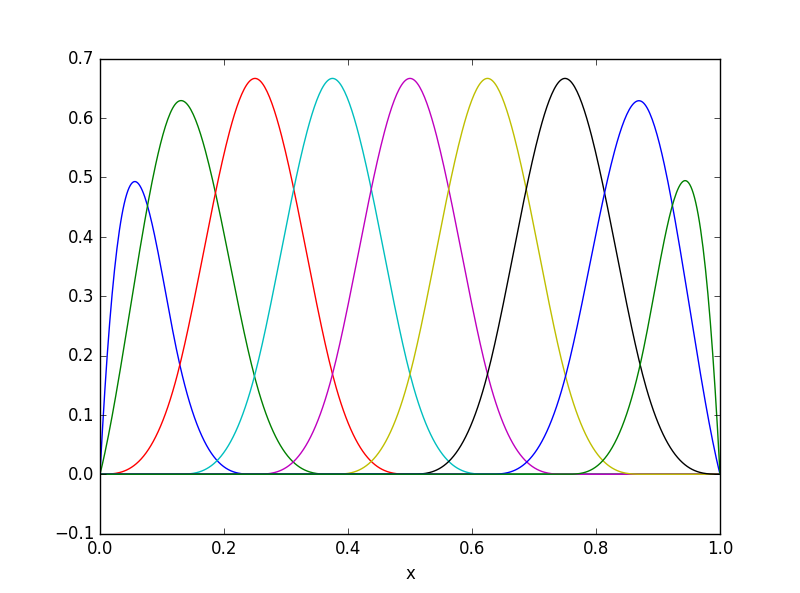}
	\end{tabular}
	\caption{Left panel: The cubic B-spline (red line) and its first derivative (blue line). Right panel: The B-basis $\{\phi_{3,3,k}(x)\}$.}
	\label{fig:Bcubic}
\end{figure}

To shown the effectiveness of the fractional spline collocation-Galerkin method we solved the fractional diffusion problem
(\ref{eq:fracdiffeq}) for two different known terms $f(t,x)$ taken from \cite{FXY11}.
In all the numerical tests we used as approximating space for the Galerkin method the (polynomial) cubic spline space. The B-splines $B_3$, its first derivatives $B_3'$ and the B-basis $\{\phi_{3,3,k}\}$ are displayed in Figure~\ref{fig:Bcubic}. We notice that since the cubic B-spline is centrally symmetric in the interval $[0,4]$, the B-basis is centrally symmetric, too. All the numerical tests were performed on a laptop using a Python environment. Each test takes a few minutes.

\subsection{Example 1}

\begin{figure}[htp]
	\centering
	\begin{tabular}{cc}
		\includegraphics[width=0.45\textwidth]{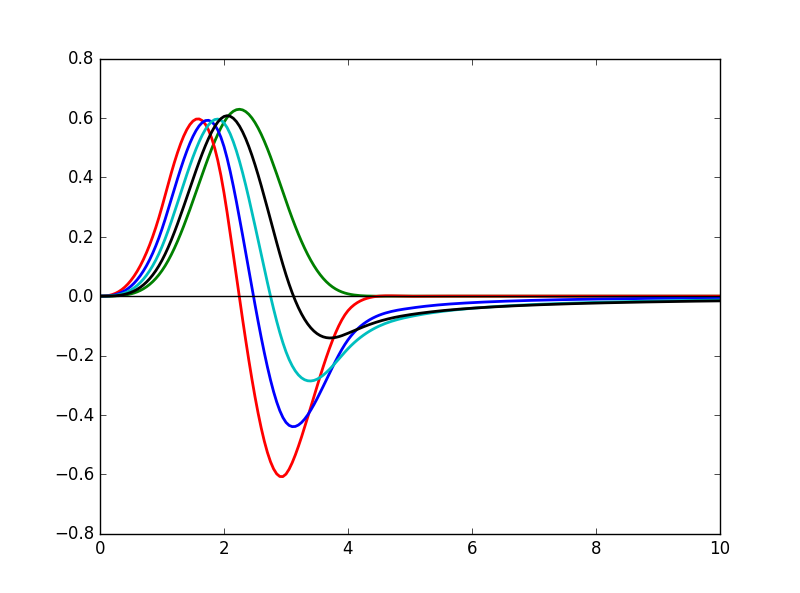} &
		\includegraphics[width=0.45\textwidth]{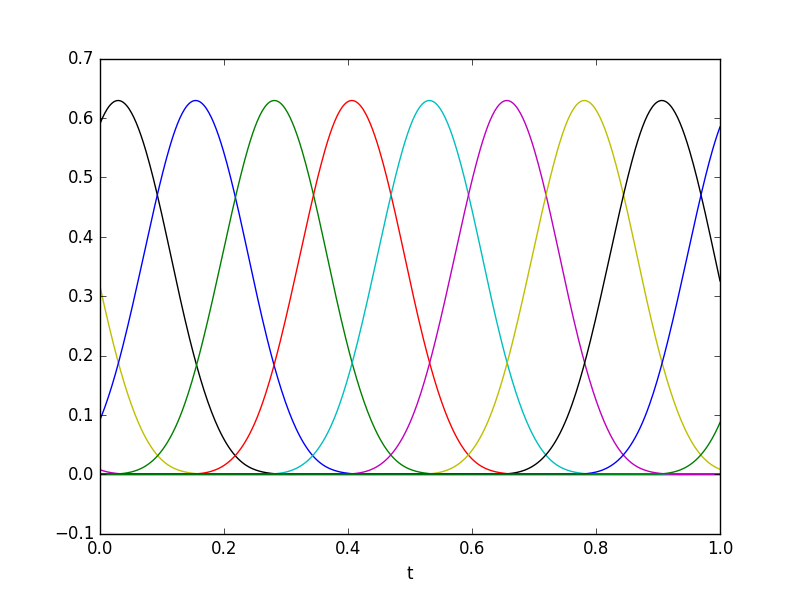}
	\end{tabular}
	\caption{Left panel: The fractional B-spline $B_{3.5}$ (green line), its first derivative (red line) and its fractional derivatives of order $\gamma =0.75$ (blue line), 0.5 (cyan line), 0.25 (black line).  Right panel: The fractional basis $\{\chi_{3.5,3,r}\}$(right).}
	\label{fig:fract_Basis}
\end{figure}

\begin{figure}[htp]
	\centering
	\begin{tabular}{cc}
		\includegraphics[width=0.45\textwidth]{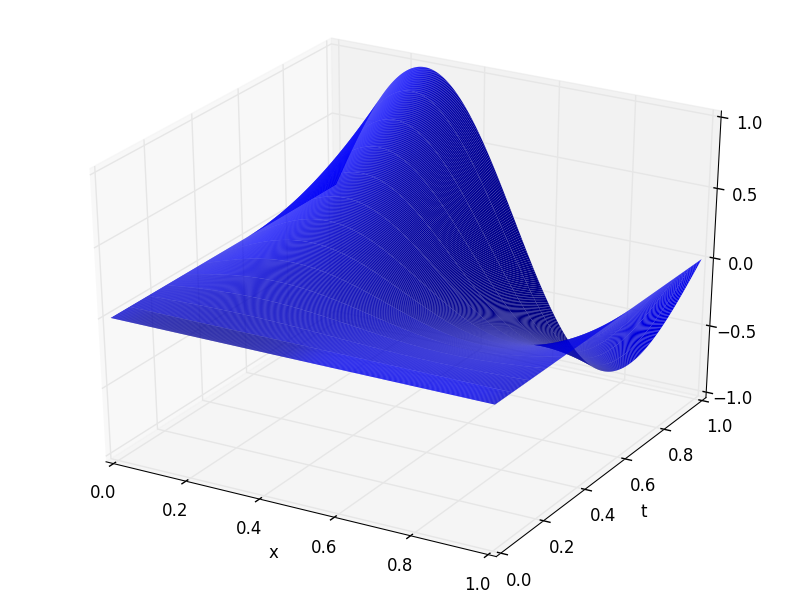} &
		\includegraphics[width=0.45\textwidth]{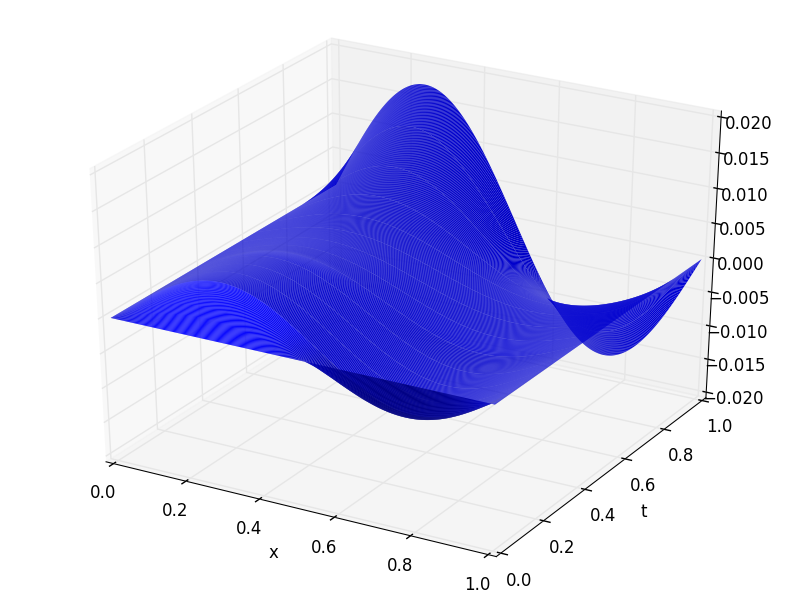} \\
	\end{tabular}
	\caption{Example 1. The numerical solution (left panel) and the error (right panel) for $j=6$ and $s=6$ when $\gamma = 0.5$.}
	\label{fig:numsol_1}
\end{figure}

In the first test we solved the time-fractional diffusion equation (\ref{eq:fracdiffeq}) in the case when  
$$f(t,x)=\frac{2}{\Gamma(3-\gamma)}\,t^{2-\gamma}\, \sin(2\pi x)+4\pi^2\,t^2\, \sin(2\pi x)\,.$$ 
The exact solution is 
$$u(t,x)=t^2\,\sin(2\pi x).$$
We used the fractional B-spline $B_{3.5}$ as approximating function for the collocation method and solved the problem for $\gamma = 1, 0.75, 0.5, 0.25$. The  fractional B-spline $B_{3.5}$, its first derivative and its fractional derivatives are shown in Figure~\ref{fig:fract_Basis} along with the fractional basis $\{\chi_{3.5,3,r}\}$. The numerical solution $u_{s,j}(t,x)$ and the error $e_{s,j}(t,x) = u(t,x)-u_{s,j}(t,x)$ for $s=6$ and $j=6$ are displayed in Figure~\ref{fig:numsol_1} for $\gamma = 0.5$. In all the numerical tests we set $q = s+1$.

\begin{figure}[htp]
	\centering
	\begin{tabular}{cc}
		\includegraphics[width=0.45\textwidth]{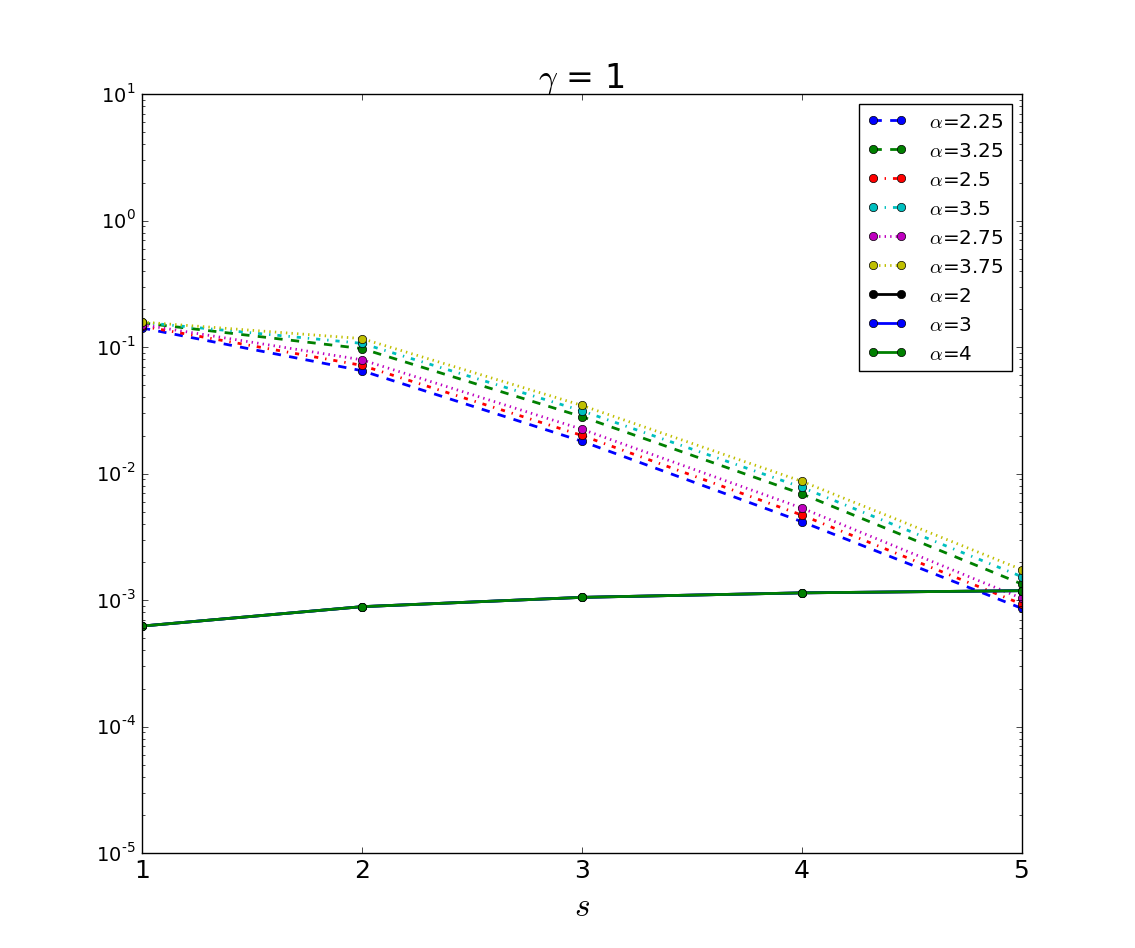} &
		\includegraphics[width=0.45\textwidth]{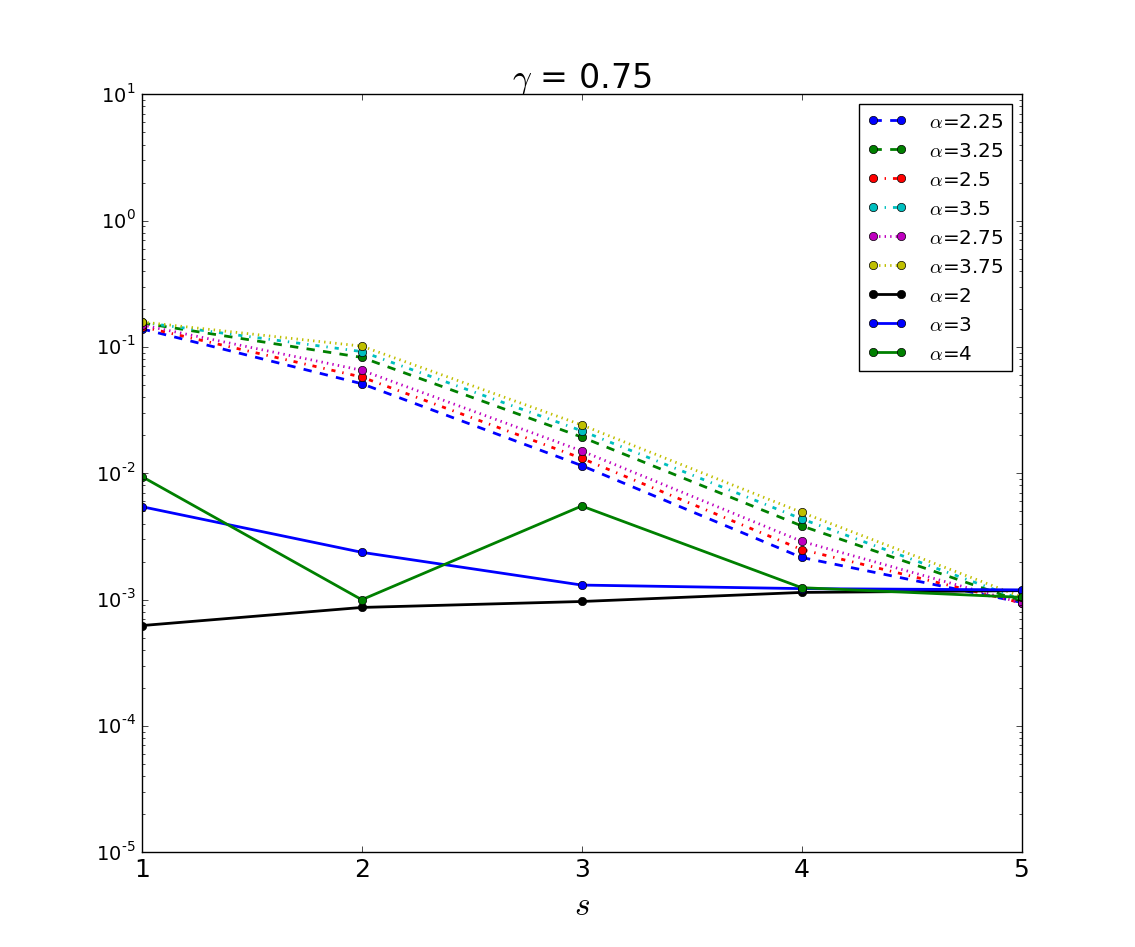}
		\\ 
		\includegraphics[width=0.45\textwidth]{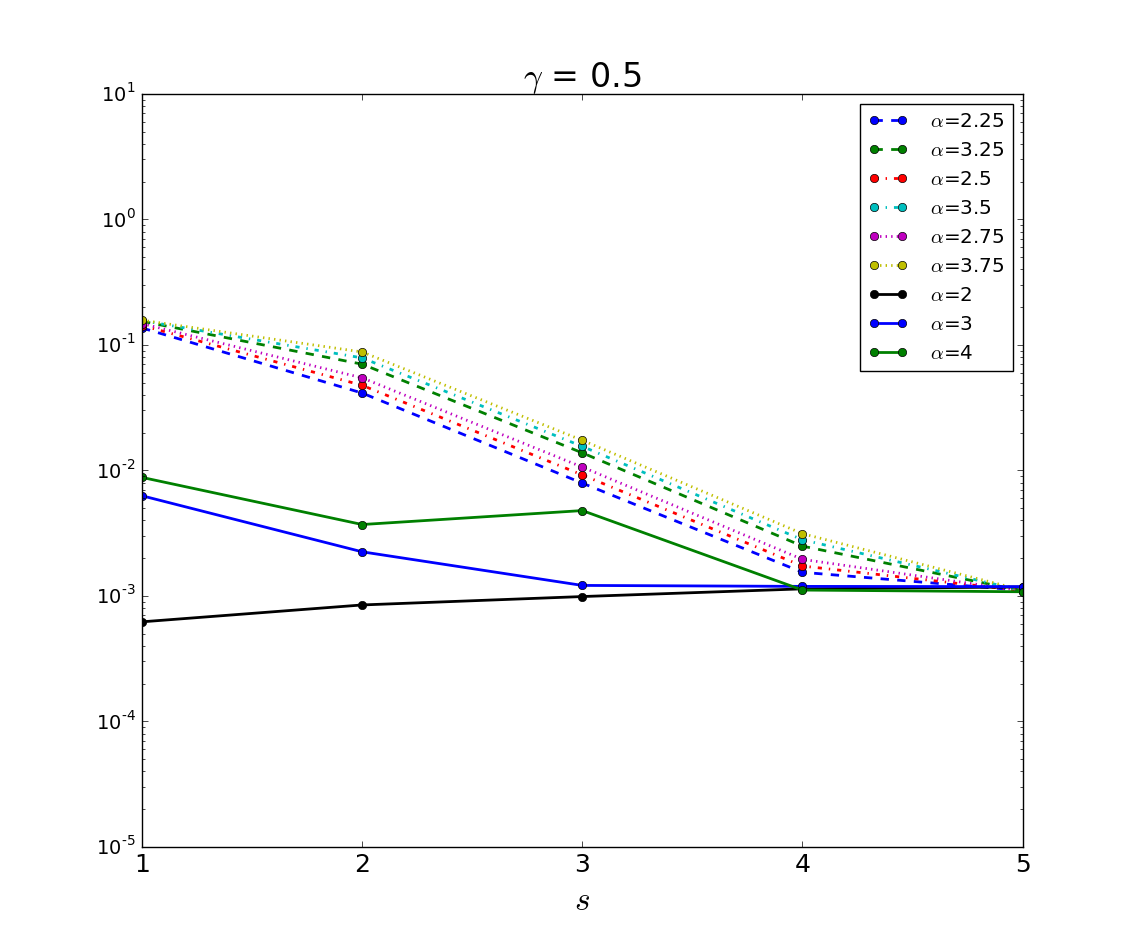} &
		\includegraphics[width=0.45\textwidth]{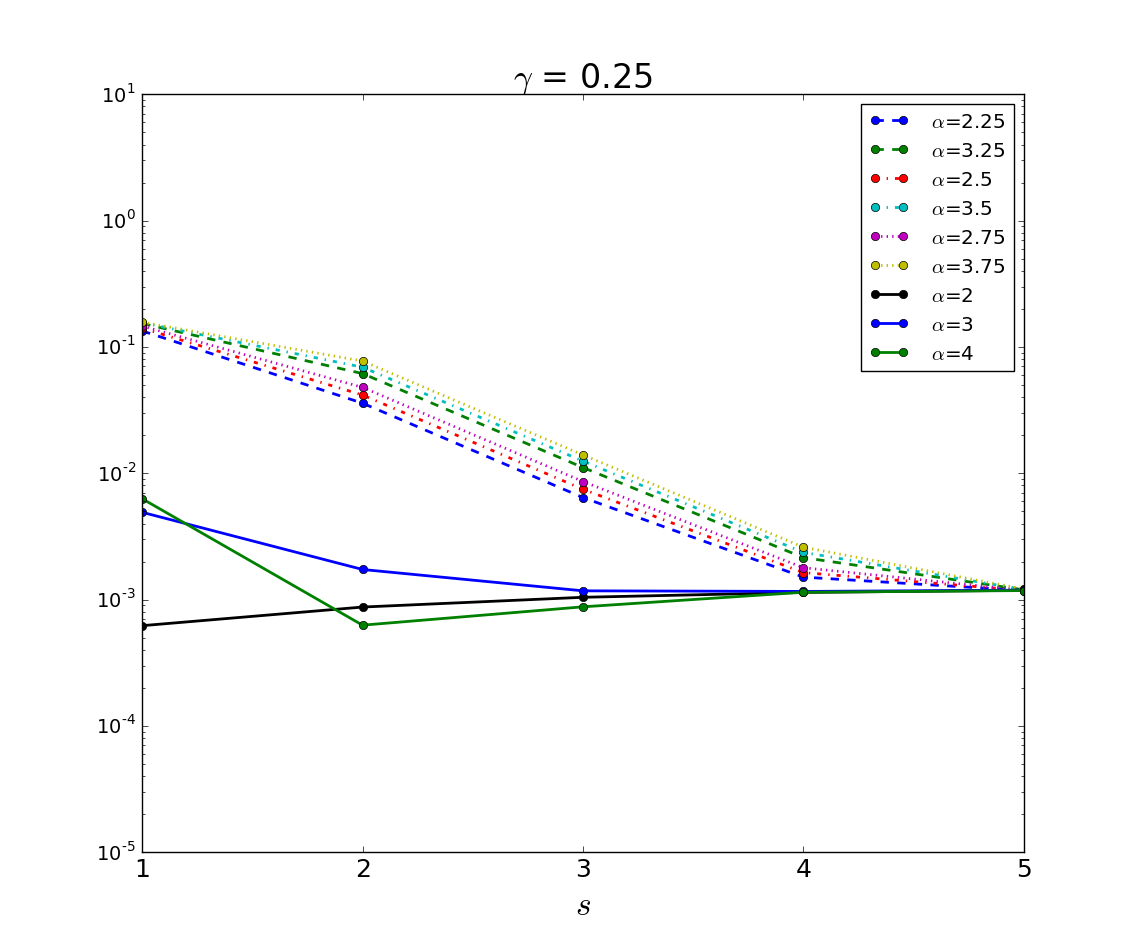}
	\end{tabular}
	\caption{Example 1: The $L_2$-norm of the error as a function of $q$ for different values of $\gamma$. Each line corresponds to a spline of different degree: solid lines correspond to the polynomial splines; non solid lines correspond to fractional splines.}
	\label{fig:L2_error_1}
\end{figure}

We analyze the behavior of the error as the degree of the fractional B-spline $B_\beta$ increases. Figure~\ref{fig:L2_error_1} shows the $L_2$-norm of the error as a function of $s$ for $\beta$ ranging from 2 to 4; the four panels in the figure refer to different values of the order of the fractional derivative. For these tests we set $j=5$. The figure shows that for $s \le 4$ the error provided by the polynomial spline approximations is lower than the error provided by the fractional spline approximations. Nevertheless, in this latter case the error decreases reaching the same value, or even a lower one, of the polynomial spline error when $s=5$. We notice that for $\gamma=1$ the errors provided by the polynomial spline approximations of different degrees have approximatively the same values while the error provided by the polynomial spline of degree 2 is lower in case of fractional derivatives. In fact, it is well-known that fractional derivatives are better approximated by less smooth functions \cite{Po99}.

Then, we analyze the convergence of the method for increasing values of $j$ and $s$. 
Table~\ref{tab:conv_js_fract_1} reports the $L_2$-norm of the error for different values of $j$ and $s$ when using the fractional B-spline $B_{3.5}$  and $\gamma = 0.5$. The number of degrees-of-freedom is also reported. The table shows that the error decreases when $j$ increases and $s$ is held fix. 
We notice that the error decreases very slightly when $j$ is held fix and $s$ increases since for these values of $s$ we reached the accuracy level we can expect for that value of $j$  (cf. Figures~\ref{fig:L2_error_1}). The higher values of the error for $s=7$ and $j=5,6$ are due to the numerical instabilities of the basis $\{\chi_{3.5,s,r}\}$ which result in a high condition number of the discretization matrix. The error has a similar behavior  even in the case when we used the cubic B-spline space as approximating space for the collocation method (cf. Table~\ref{tab:conv_js_cubic_1}).

\begin{table}[htp]
	\caption{Example 1: The $L_2$-norm of the error for increasing values of $s$ and $j$ when using the fractional B-spline of degree $\beta=3.5$. The numbers in parenthesis are the degrees-of-freedom. Here, $\gamma =0.5$.}
	{\small \begin{tabular}{c|c|c|c|c}
		\backslashbox{$s$}{$j$} & 3 & 4 & 5 & 6 \\
		\hline \hline
		$\sharp V_j^{(\alpha)}([0,1])$ & 9 & 17 & 33 & 65 \\
		\hline \hline
		5 & 0.02037 (369)  & 0.00449 (697)  & 0.00101 (1353) & 0.00025 (2665)\\
		6 & 0.02067 (657)  & 0.00417 (1241) & 0.00093 (2409) & 0.00024 (4745)\\
		7 & 0.01946 (1233) & 0.00381 (2329) & 0.00115 (4521) & 0.00117 (8905)
		\\
		\hline
	\end{tabular}}
	\label{tab:conv_js_fract_1}
\end{table}

\begin{table}[htp]
	\caption{Example 1: The $L_2$-norm of the error for increasing values of $s$ and $j$ when using the cubic B-spline. The numbers in parenthesis are the degrees-of-freedom.  Here, $\gamma =0.5$.}
	{\small \begin{tabular}{c|c|c|c|c}
		\backslashbox{$s$}{$j$} & 3 & 4 & 5 & 6 \\
		\hline \hline
		$\sharp V_j^{(\alpha)}([0,1])$ & 9 & 17 & 33 & 65 \\
		\hline \hline
		5 & 0.02121 (315)  & 0.00452 (595)  & 0.00104 (1155) & 0.00025 (2275)\\
		6 & 0.02109 (603)  & 0.00443 (1139) & 0.00097 (2211) & 0.00023 (4355)\\
		7 & 0.02037 (1179) & 0.00399 (2227) & 0.00115 (4323) & 0.00115 (8515)
		\\
		\hline
	\end{tabular}}
	\label{tab:conv_js_cubic_1}
\end{table}

\subsection{Example 2}

\begin{figure}[htp]
	\centering
	\begin{tabular}{cc}
		\includegraphics[width=0.45\textwidth]{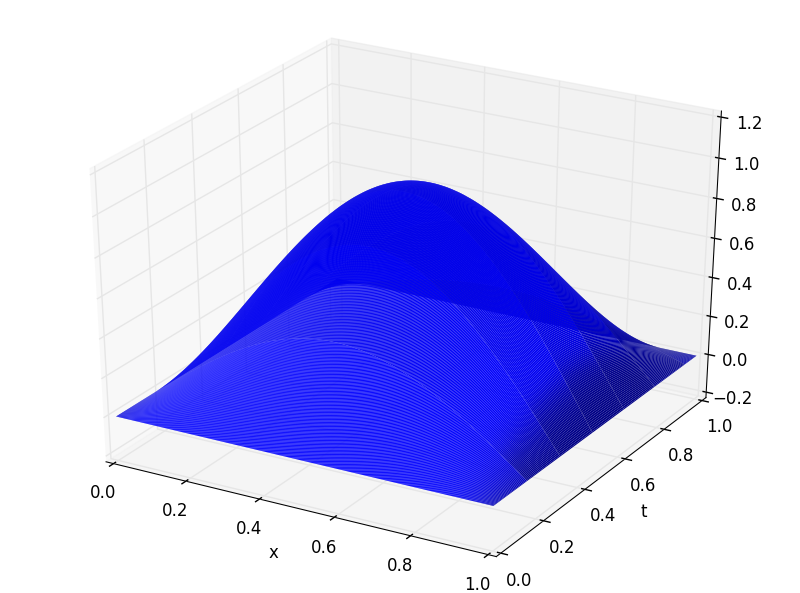} &
		\includegraphics[width=0.45\textwidth]{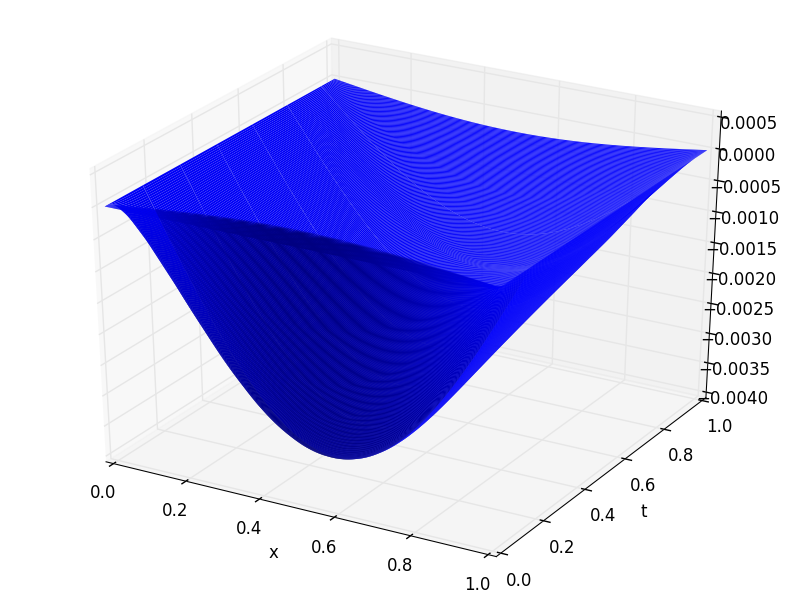} \\
	\end{tabular}
	\caption{Example 2: The numerical solution (left panel) and the error (right panel) when $j=6$ and $s=5$. }
	\label{fig:numsol_2}
\end{figure}

\begin{figure}[htp]
	\centering
	\begin{tabular}{cc}
		\includegraphics[width=0.45\textwidth]{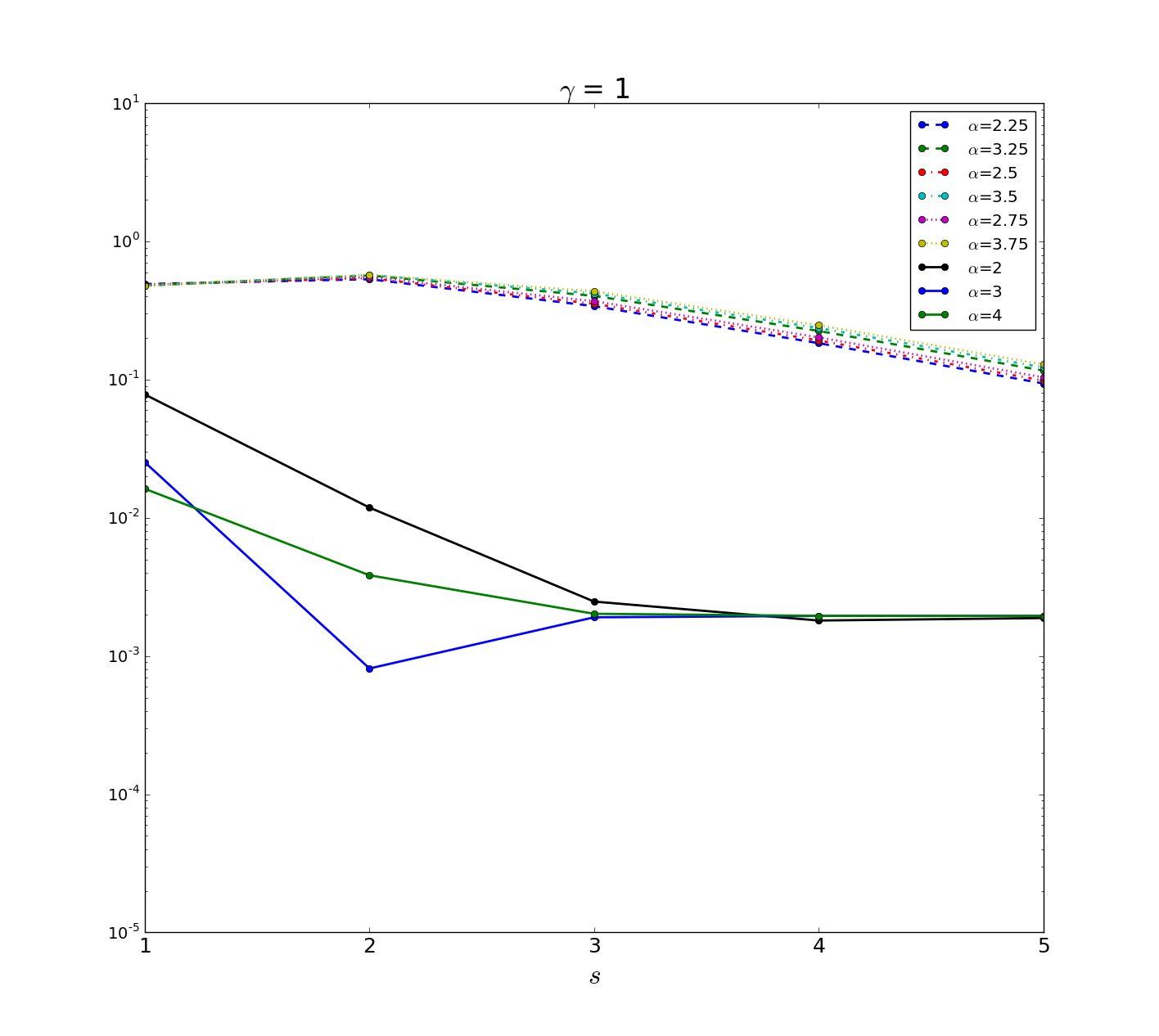} &
		\includegraphics[width=0.45\textwidth]{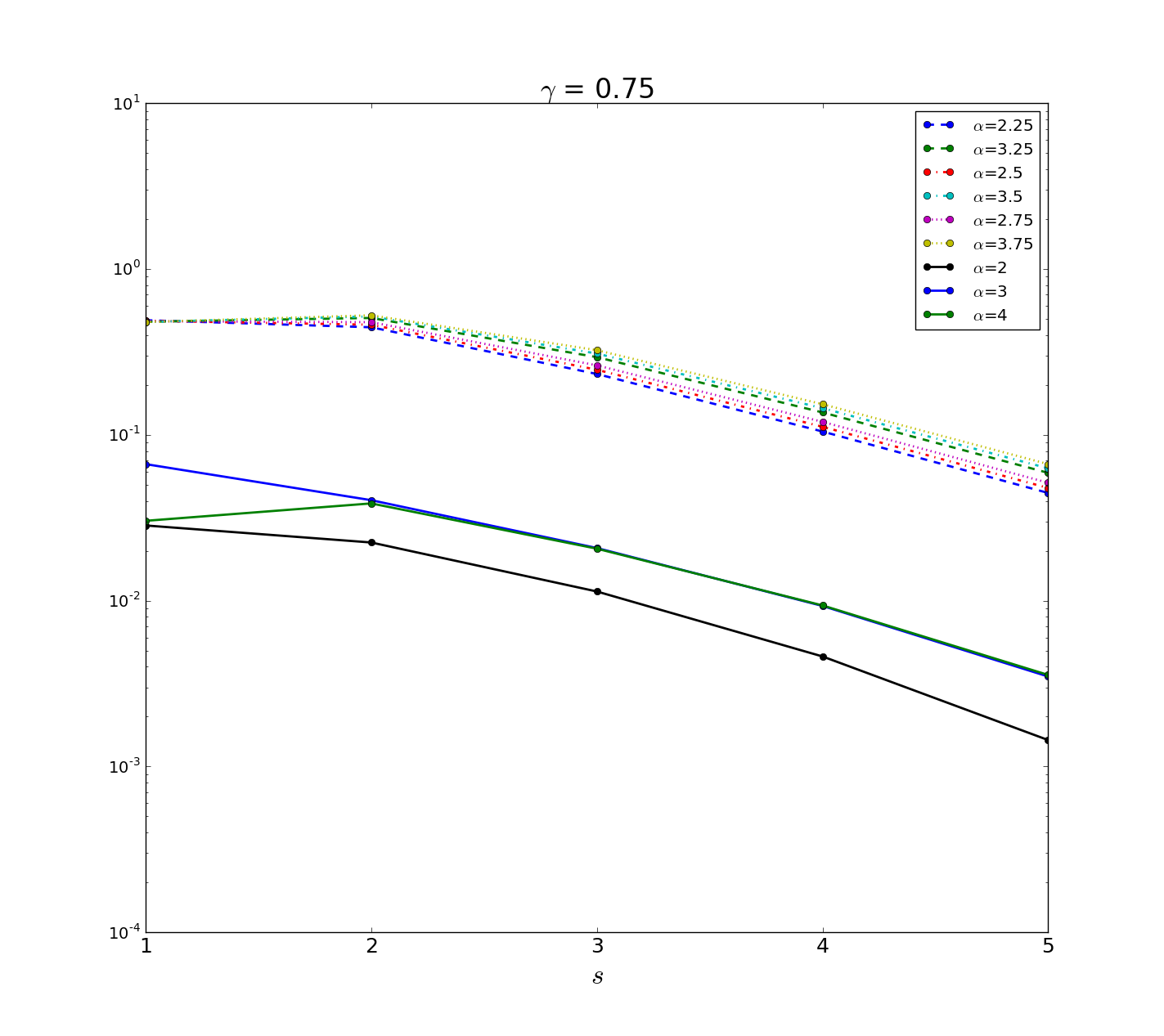}
		\\ 
		\includegraphics[width=0.45\textwidth]{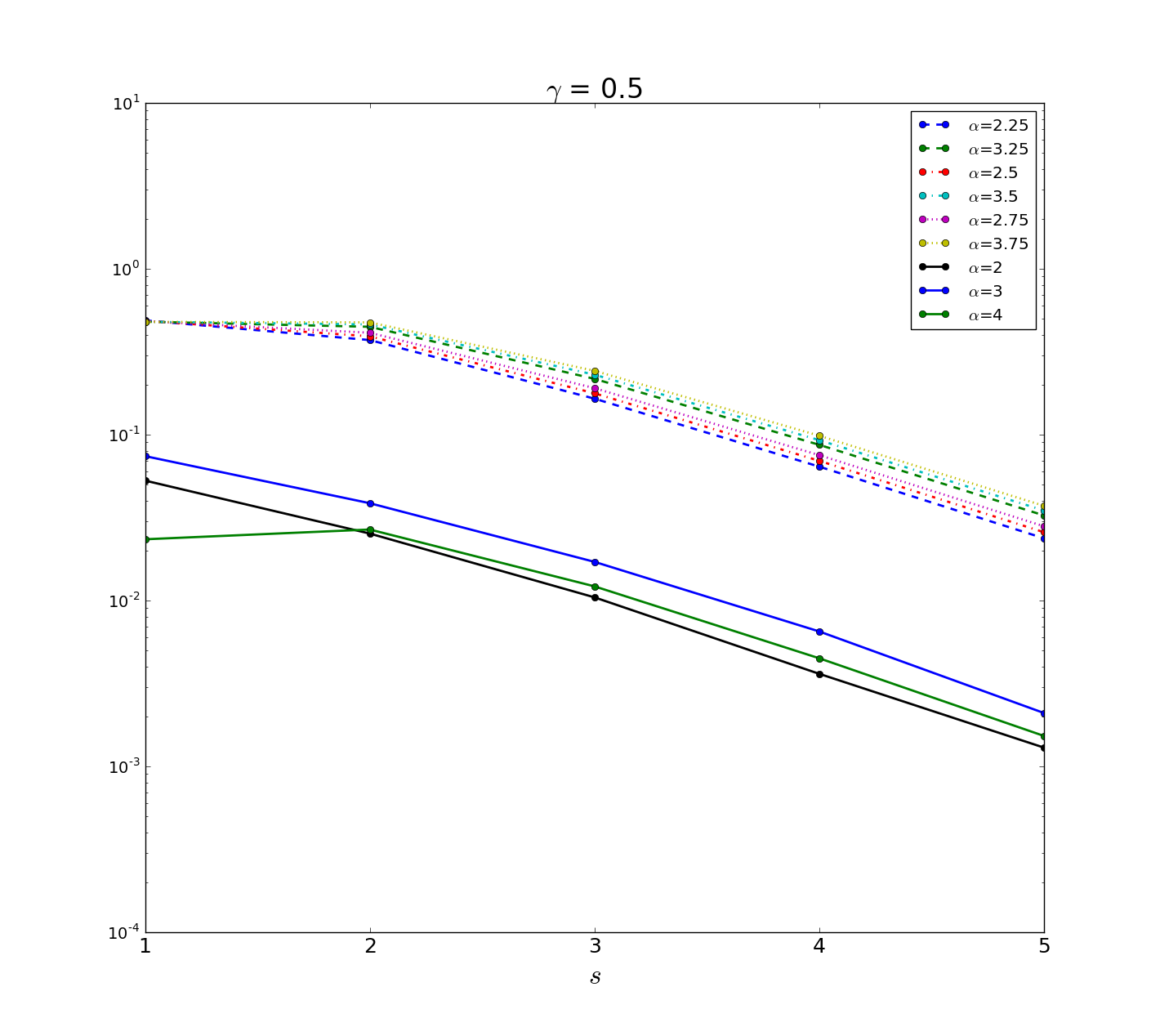} &
		\includegraphics[width=0.45\textwidth]{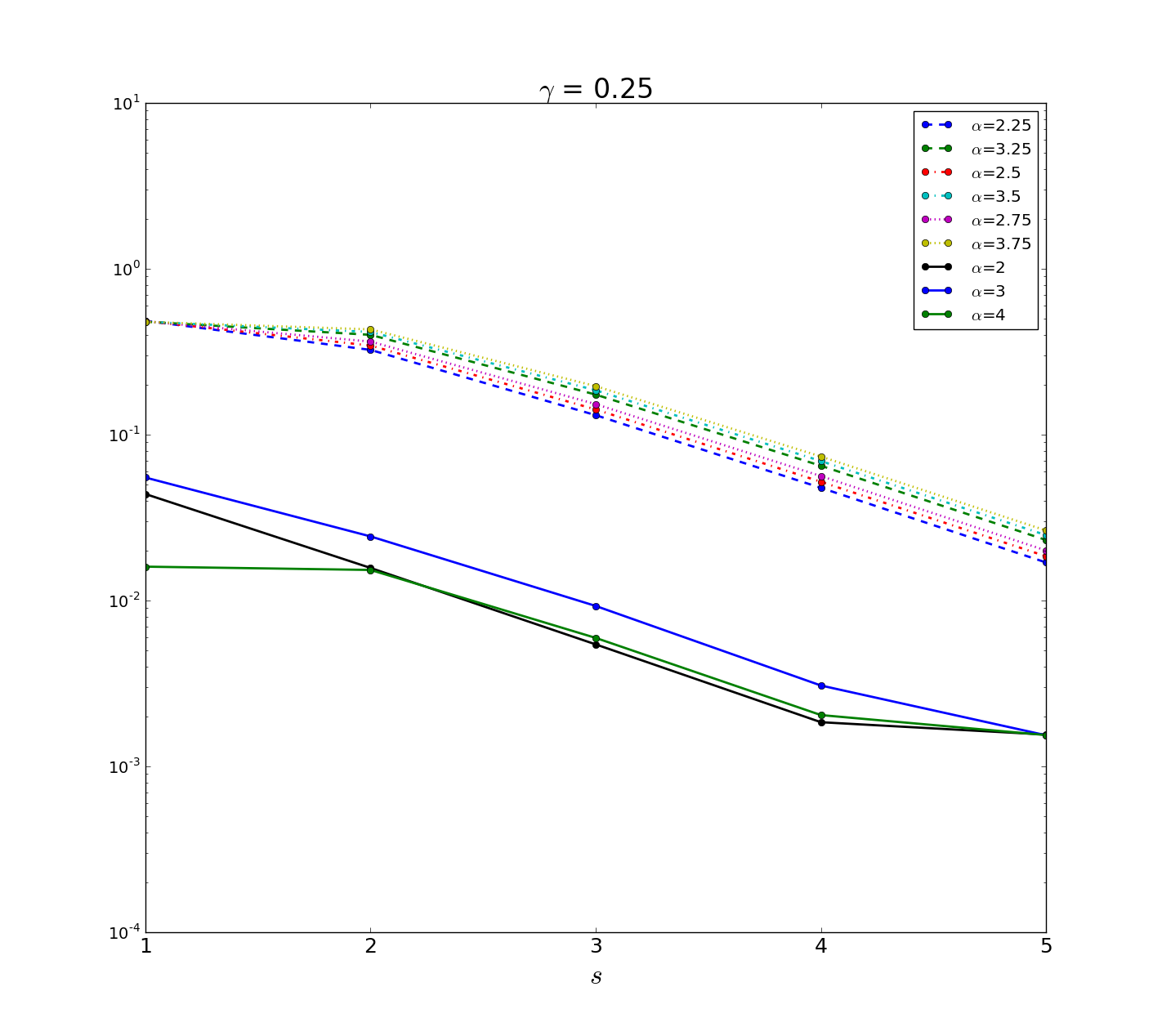}
	\end{tabular}
	\caption{Example 2: The $L_2$-norm of the error as a function of $s$ for different values of $\gamma$. Each line corresponds to a spline of different degree: solid lines correspond to the polynomial splines; non solid lines correspond to fractional splines.}
	\label{fig:L2_error_2}
\end{figure}

In the second test we solved the time-fractional diffusion equation (\ref{eq:fracdiffeq}) in the case when  
$$\begin{array}{lcl}
f(t,x) & = &  \displaystyle \frac{\pi t^{1-\gamma}}{2\Gamma(2-\gamma)} \left( \, {_1F_1}(1,2-\gamma,i\pi\,t) + \,{_1F_1}(1,2-\gamma,-i\pi\,t) \right) \, \sin(\pi\,x)  \\ \\
& + & \pi^2 \, \sin(\pi\,t)  \, \sin(\pi\,x)\,,
\end{array}
$$
where $_1F_1(\alpha,\beta,z)$ is the Kummer's confluent hypergeometric function  defined as
$$
_1F_1(\alpha,\beta, z) =
\frac {\Gamma(\beta)}{\Gamma(\alpha)} \, \sum_{k\in \NN} \, \frac {\Gamma(\alpha+k)}{\Gamma(\beta+k)\, k!} \, z^k\,, \qquad \alpha \in \RR\,, \quad -\beta \notin \NN_0\,,
$$
where $\NN_0 = \NN \backslash \{0\}$ (cf. \cite[Chapter 13]{AS65}).
In this case the exact solution is $$u(t,x)=\sin(\pi t)\,\sin(\pi x).$$
We performed the same set of numerical tests as in Example 1. 
The numerical solution $u_{s,j}(t,x)$ and the error $e_{s,j}(t,x)$ for $s=5$ and $j=6$ are displayed in Figure~\ref{fig:numsol_2} in the case when $\gamma = 0.5$.
Figure~\ref{fig:L2_error_2} shows the $L_2$-norm of the error as a function of $s$ for $\beta$ ranging from 2 to 4 and $j=5$; the four panels in the figure refer to different values of the order of the fractional derivative. 
Tables~\ref{tab:conv_js_fract_2}-\ref{tab:conv_js_cubic_2} report the $L_2$-norm of the error for different values of $j$ and $s$ and $\beta = 3.5, 3$, respectively. The number of degrees-of-freedom is also reported. 
\\
Figure~\ref{fig:L2_error_2} shows that value of the error is higher than in the previous example but it decreases as $s$ increases showing a very similar behavior as that one in Example 1. The values of the error in Tables~\ref{tab:conv_js_fract_2}-\ref{tab:conv_js_cubic_2} are approximatively the same as in Tables~\ref{tab:conv_js_fract_1}-\ref{tab:conv_js_cubic_1}.

\begin{table}[htp]
	\caption{Example 2: The $L_2$-norm of the error for increasing values of $s$ and $j$  for the fractional B-spline of degree $\beta=3.5$. The numbers in parenthesis are the degrees-of-freedom. Here, $\gamma=0.5$.}
	{\small \begin{tabular}{c|c|c|c|c}
			\backslashbox{$s$}{$j$} & 3 & 4 & 5 & 6 \\
			\hline \hline
			$\sharp V_j^{(\alpha)}([0,1])$ & 9 & 17 & 33 & 65 \\
			\hline \hline
			5 & 0.01938 (369)  & 0.00429 (697)  & 0.00111 (1353) & 0.00042 (2665) \\
			6 & 0.01809 (657)  & 0.00555 (1241) & 0.00507 (2409) & 0.00523 (4745) \\
			7 & 0.01811 (1233) & 0.01691 (2329) & 0.01822 (4521) & 0.01858 (8905)
			\\
			\hline
	\end{tabular}}
	\label{tab:conv_js_fract_2}
\end{table}

\begin{table}[htp]
	\caption{Example 2: The $L_2$-norm of the error for increasing values of $s$ and $j$ for the cubic B-spline. The numbers in parenthesis are the degrees-of-freedom. Here, $\gamma=0.5$.}
	{\small \begin{tabular}{c|c|c|c|c}
			\backslashbox{$s$}{$j$} & 3 & 4 & 5 & 6 \\
			\hline \hline
			$\sharp V_j^{(\alpha)}([0,1])$ & 9 & 17 & 33 & 65 \\
			\hline \hline
			5 & 0.01909 (315)  & 0.00404 (595)  & 0.00102 (1155) & 0.00063 (2275) \\
			6 & 0.01810 (603)  & 0.00546 (1139) & 0.00495 (2211) & 0.00511 (4355) \\
			7 & 0.01805 (1179) & 0.01671 (2227) & 0.01801 (4323) & 0.01838 (8515)
			\\
			\hline
	\end{tabular}}
	\label{tab:conv_js_cubic_2}
\end{table}

\section{Conclusion}
\label{sec:concl}

We proposed a fractional spline collocation-Galerkin method to solve the time-fractional  diffusion equation. The novelty of the method is in the use of fractional spline spaces as approximating spaces so that the fractional derivative of the approximating function can be evaluated easily by 
an explicit differentiation rule that involves the generalized finite difference operator. The numerical tests show that the method has a good accuracy
so that it can be effectively used to solve fractional differential problems. The numerical instabilities arising in the fractional basis when $s$ increases can be reduced following the approach in \cite{GPP04} that allows us to construct stable basis on the interval. Moreover, the ill-conditioning of the linear system (\ref{colllinearsys}) can be reduced using iterative methods in Krylov spaces, such as the method proposed in \cite{CPSV17}.
Finally, we notice that following the procedure given in \cite{GPP04},  fractional wavelet bases on finite interval can be constructed so that the proposed method can be generalized to fractional wavelet approximating spaces.



\end{document}